2000]20F67

# THE GEOMETRY OF RELATIVE CAYLEY GRAPHS FOR SUBGROUPS OF HYPERBOLIC GROUPS

ILYA KAPOVICH

Abstract. We show that if $H$ is a quasiconvex subgroup of a hyperbolic group $G$ then the relative Cayley graph $Y$ (also known as the Schreier coset graph) for $G/H$ is Gromov-hyperbolic. We also observe that in this situation if $G$ is torsion-free and non-elementary and $H$ has infinite index in $G$ then the simple random walk on $Y$ is transient.

## 1. Introduction

We will call a map $\pi : A \to G$ a *marked finite generating set of* $G$ if $A$ is a finite alphabet disjoint from $G$ and the set $\pi(A)$ generates $G$. (Thus we allow $\pi$ to be non-injective and also allow $1 \in \pi(G)$.) By abuse of notation we will often suppress the map $\pi$ from consideration and talk about $A$ being a marked finite generating set or just a finite generating set of $G$. We recall the explicit construction of a relative coset Cayley graph, since it is important for our considerations.

**Definition 1.1.** Let $G$ be a group and let $\pi : A \to G$ be a marked finite generating set of $G$. Let $H \le G$ be a subgroup of $G$ (not necessarily normal). The *relative Cayley graph* (or the *Schreier coset graph*) $\Gamma(G/H, A)$ for $G$ relative $H$ with respect to $A$ is an oriented labeled graph defined as follows:

1. The vertices of $\Gamma(G/H, A)$ are precisely the cosets $G/H = \{Hg \mid g \in G\}$.
2. The set of positively oriented edges of $\Gamma(G/H, A)$ is in one-to-one correspondence with the set $G/H \times A$. For each pair $(Hg, a) \in G/H \times A$ there is a positively oriented edge in $\Gamma(G/H, A)$ from $Hg$ to $Hg\pi(a)$ labeled by the letter $a$.

Thus the label of every path in $\Gamma(G/H, A)$ is a word in the alphabet $A \cup A^{-1}$. The graph $\Gamma(G/H, A)$ is connected since $\pi(A)$ generates $G$. Moreover, $\Gamma(G/H, A)$ comes equipped with the natural simplicial metric $d$ obtained by giving every edge length one.

The relative Cayley graph $\Gamma(G/H, A)$ can be identified with the 1-skeleton of the covering space corresponding to $H$ of the presentation complex of $G$ on the generating set $A$. Note that if $H$ is normal in $G$ and $G_1 = G/H$ is the quotient group, then $\Gamma(G/H, A)$ is exactly the Cayley graph of the group $G_1$ with respect to $A$. In particular, if $H = 1$ then $\Gamma(H/1, A)$ is the standard Cayley graph of $G$ with respect to $A$, denoted $\Gamma(G, A)$.

Recall that a geodesic metric space $(X, d)$ is said to be $\delta$-*hyperbolic* if every geodesic triangle in $X$ is $\delta$-*thin* that is each side of the triangle is contained in the closed $\delta$-neighborhood of the union of the other sides. A geodesic space $(X, d)$ is *hyperbolic* (or *Gromov-hyperbolic*) if it is $\delta$-hyperbolic for some $\delta \ge 0$. According to M.Gromov [30], a finitely generated group $G$ is said to be *hyperbolic* (or *word-hyperbolic*) if for any marked finite generating set $A$ of $G$ the Cayley graph $\Gamma(G, A)$ is a hyperbolic metric space. Recall that a subgroup $H$ of a word-hyperbolic group $G$ is said to be *quasiconvex* if for any finite generating set $A$ of $G$ there is $C > 0$ such that every geodesic in $\Gamma(G, A)$ with both endpoints in $H$ is contained in the $C$-neighborhood of $H$. Quasiconvex subgroups play a very important role in the theory of hyperbolic groups, as they enjoy some particularly good properties. Thus a quasiconvex subgroup $H$ of a hyperbolic group $G$ is itself finitely generated, finitely presentable and word-hyperbolic. Moreover, $H$ has solvable membership problem with in $G$ and the subgroup respect $H$ is rational with respect to any automatic structure on $G$. The intersection of any finite family of quasiconvex subgroups of $G$ is again quasiconvex (and hence finitely generated). It is also known that a finitely generated subgroup $H$ of a hyperbolic group $G$ is quasiconvex if and only if $H$ is quasi-isometrically embedded in $G$. Quasiconvex subgroups are crucial in the Combination Theorem for hyperbolic groups [4, 38, 24], as well as its various applications (see for example [34, 42, 21, 35, 36, 55]). Quasiconvexity for subgroups of hyperbolic groups is closely related to geometric finiteness for Kleinian groups [54, 45, 37, 55]. Relative Cayley graphs are important for understanding the pro-finite topology on $G$ and investigating the separability properties for quasiconvex subgroups of hyperbolic groups, as demonstrated by the work of R.Gitik [26, 27, 28]. If $H \le G$ then the number of relative ends $e(G,H)$ is defined as the number of ends of the coset graph $\Gamma(G/H, A)$. Therefore the study of connectivity at infinity of coset graphs also arises naturally in various generalizations of Stallings' theory about ends of groups [48, 19, 12, 50]. Studying the situation, when $H$ is a virtually cyclic subgroup (and hence quasiconvex) of a hyperbolic group $G$ and $e(G,H) > 1$, is crucial in the theory of JSJ-decomposition for word-hyperbolic groups developed by Z.Sela [52] and later by B.Bowditch [6] (see also [47, 12, 18, 51] for various generalizations of the theory of JSJ-splittings). The general case when $H$ is a quasiconvex subgroup of a hyperbolic group $G$ and $e(G,H) > 1$ was addressed by M.Sageev in [49], where it leads to actions of $G$ on finite-dimensional non-positively curved cubings. Understanding the geometry of the coset graphs for quasiconvex subgroups is particularly important for computational purposes, such as performing the Todd-Coxeter process, solving the uniform membership problem with respect to quasiconvex subgroups and finding the quasiconvexity constant algorithmically [33, 32, 14, 15]. Moreover, looking at relative Cayley graphs (particularly for the case of a quasiconvex subgroup in a hyperbolic or automatic group) is central for the study of the topologically motivated notion of a *tame subgroup*, which was pursued by M.Mihalik [39, 40] and R.Gitik [28].

Our main result is the following:

**Theorem 1.2.** *Let $G$ be a word-hyperbolic group with a marked finite generating set $A$. Let $H \le G$ be a quasiconvex subgroup. Then the relative Cayley graph $Y = \Gamma(G/H, A)$ is a hyperbolic metric space.*

The statement of Theorem 1.2 appears as a claim without proof in Section 5.3, page 139 of M.Gromov's book [30]. Theorem 1.2 was also independently obtained by Robert Foord [17] in his PhD thesis at Warwick (which appears unlikely to be published according to the information received by the author from Foord's advisor Derek Holt).

It must be stressed that hyperbolicity of the coset graph $\Gamma(G/H, A)$ has nothing to do with relative hyperbolicity of $H$ in $G$, as defined by either M.Gromov, B.Farb or B.Bowditch [31, 16, 5, 56]. Relative hyperbolicity in either sense roughly speaking deals with collapsing in $\Gamma(G, A)$ into single points all pairs $g, gh$ (where $g \in G, h \in H$). For example, if $G = F(a, b)$ and $H = \langle a \rangle$, then the length of the path $p = aba^{100}b^{-1}a^{100}$ essentially becomes 2. To construct the relative Cayley graph $\Gamma(G/H, A)$ from $\Gamma(G, A)$ one has to collapse the *right cosets $Hg$ only.* Thus for the above example the distance in $\Gamma(G/H, A)$ between $H1$ and $Hp$ is equal to 202. The graph $\Gamma(G/H, A)$ generally does not admit a $G$-action, which is very much unlike the relatively hyperbolic situation. We should also stress that B.Bowditch's general technique of collapsing a family of separated uniformly quasiconvex subsets in a hyperbolic space to obtain a new hyperbolic space [5] is not applicable for our purposes, since the sets $Hg$ are not uniformly quasiconvex in $\Gamma(G, A)$. If $H$ is quasiconvex in a hyperbolic group $H$, then $G$ is known to be relatively hyperbolic with respect to $H$ in the senses of both B.Farb and B.Bowditch (see [16, 20, 5]).

It is worth noting that the statement of Theorem 1.2 does not hold for arbitrary finitely generated subgroups of hyperbolic groups. For example, by a remarkable result of E.Rips [46] for any finitely presented group $Q$ there exists a short exact sequence

$$1 \to K \to G \to Q \to 1$$

where $G$ is torsion-free non-elementary word-hyperbolic and where $K$ is two-generated and non-elementary. Thus $G/K = Q$ and if $Q$ is not hyperbolic to begin with, then the relative Cayley graph $\Gamma(G/K, A)$ is not hyperbolic either. On the other hand it is possible that $H \le G$ is not quasiconvex but $\Gamma(G/H, A)$ is Gromov-hyperbolic nonetheless. For example, this happens if $H = K$ and $Q = \mathbb{Z}$ in the construction above. Then $\Gamma(G/K, A)$ is hyperbolic since $\mathbb{Z}$ is word-hyperbolic. However, $K \le G$ is not quasiconvex, since $K$ is infinite and has infinite index in its normalizer. Similar examples can be constructed using mapping tori of automorphisms of free groups and surface groups.

Note that by construction $\Gamma(G/H, A)$ is a $2k$-regular graph where $k$ is the number of elements in $A$. We also obtain the following useful fact.

**Theorem 1.3.** *Let $G$ be a torsion-free non-elementary word-hyperbolic group with a marked finite generating set $A$. Let $H \le G$ be a quasiconvex subgroup and let $Y = \Gamma(G/H, A)$ be the relative Cayley graph.*

*Then the simple random walk on $Y$ is transient if and only if $[G : H] = \infty$.*

The statement of Theorem 1.3 no longer holds for arbitrary finitely generated subgroups of torsion-free hyperbolic groups, even if those subgroups are themselves hyperbolic. For example, if $F$ is a finitely generated free group, $\phi$ is an automorphism of $F$ without periodic conjugacy classes, then the mapping torus group $G = \langle F, t \,|\, t^{-1}ft = \phi(t), f \in F\rangle$ is word-hyperbolic (see [4, 8]). In this case $F$ is a normal subgroup of $G$ with $G/F = \mathbb{Z}$. It is, of course, well-known that the simple random walk on $\mathbb{Z}$ is recurrent.

The author is very grateful to Laurent Bartholdi, Philip Bowers, and Tatiana Smirnova-Nagnibeda for the many useful discussions regarding random walks, to Paul Schupp for the many helpful remarks and suggestions and to Derek Holt for alerting the author to the work of R.Foord.

## 2. Hyperbolic metric spaces

We assume that the reader is familiar with the basics of Gromov-hyperbolic spaces and word-hyperbolic groups (see [30, 10, 23, 1, 13] for details). We will briefly recall the main definitions.

If $(X, d)$ is a geodesic metric space and $x, y \in X$, we shall denote by $[x, y]$ a geodesic segment from $x$ to $y$ in $X$. Also, of $p$ is a path in $X$, we will denote the length of $p$ by $l(p)$. Two paths $\alpha$ and $\beta$ in $X$ are said to be *$K$-Hausdorff close* if each of them is contained in the closed $K$-neighborhood of the other.

Given a path $\alpha : [0, T] \to X$ we shall often identify $\alpha$ with its image $\alpha([0, T]) \subseteq X$. If $Z \subseteq X$ and $\epsilon \ge 0$, we will denote the closed $\epsilon$-neighborhood of $Z$ in $X$ by $N_\epsilon(Z)$.

**Definition 2.1** (Hyperbolic metric space)**.** Let $(X, d)$ be a geodesic metric space and let $\delta \ge 0$. The space $X$ is said to be *$\delta$-hyperbolic* if for any geodesic triangle in $X$ with sides $\alpha, \beta, \gamma$ we have

$$\alpha \subseteq N_\delta(\beta \cup \gamma),\ \ \beta \subseteq N_\delta(\alpha \cup \gamma)\ \text{ and } \gamma \subseteq N_\delta(\alpha \cup \beta)$$

that is for any $p \in \alpha$ there is $q \in \beta \cup \gamma$ such that $d(p, q) \le \delta$ (and the symmetric condition holds for any $p \in \beta$ and any $p \in \gamma$) .

A geodesic space $X$ is said to be *hyperbolic* if it is $\delta$-hyperbolic for some $\delta \ge 0$.

Suppose $\alpha$ and $\beta$ are geodesic segments from $x$ to $y$ in a geodesic metric space $(X, d)$. We will say that the geodesic bigon $\Theta = \alpha \cup \beta$ is *$\delta$-thin* if $\alpha \subseteq N_\delta(\beta)$ and $\beta \subseteq N_\delta(\alpha)$.

The following important result was obtained by P.Papasoglu [44]:

**Theorem 2.2.** *Let $(X, d)$ be a connected graph with the simplicial metric (that is a path-metric where the length of every edge is equal to one).*

*Then $X$ is hyperbolic if and only if there is some $\delta \ge 0$ such that all geodesic bigons in $X$ are $\delta$-thin.*

Notice that the above statement fails for arbitrary geodesic metric spaces. For example, in the Euclidean plane with the standard Euclidean metric all geodesic bigons are 0-thin, while the plane is certainly not hyperbolic.

**Remark 2.3.** P.Papasoglu only states Theorem 2.2 in [44] for the case when $X$ is the Cayley graph of a finitely generated group. However, it is easy to see that his proof does not use the group assumption at all and works for any connected graph with edges of length one. This fact was noticed by many researchers, for example W.Neumann and M.Shapiro [43].

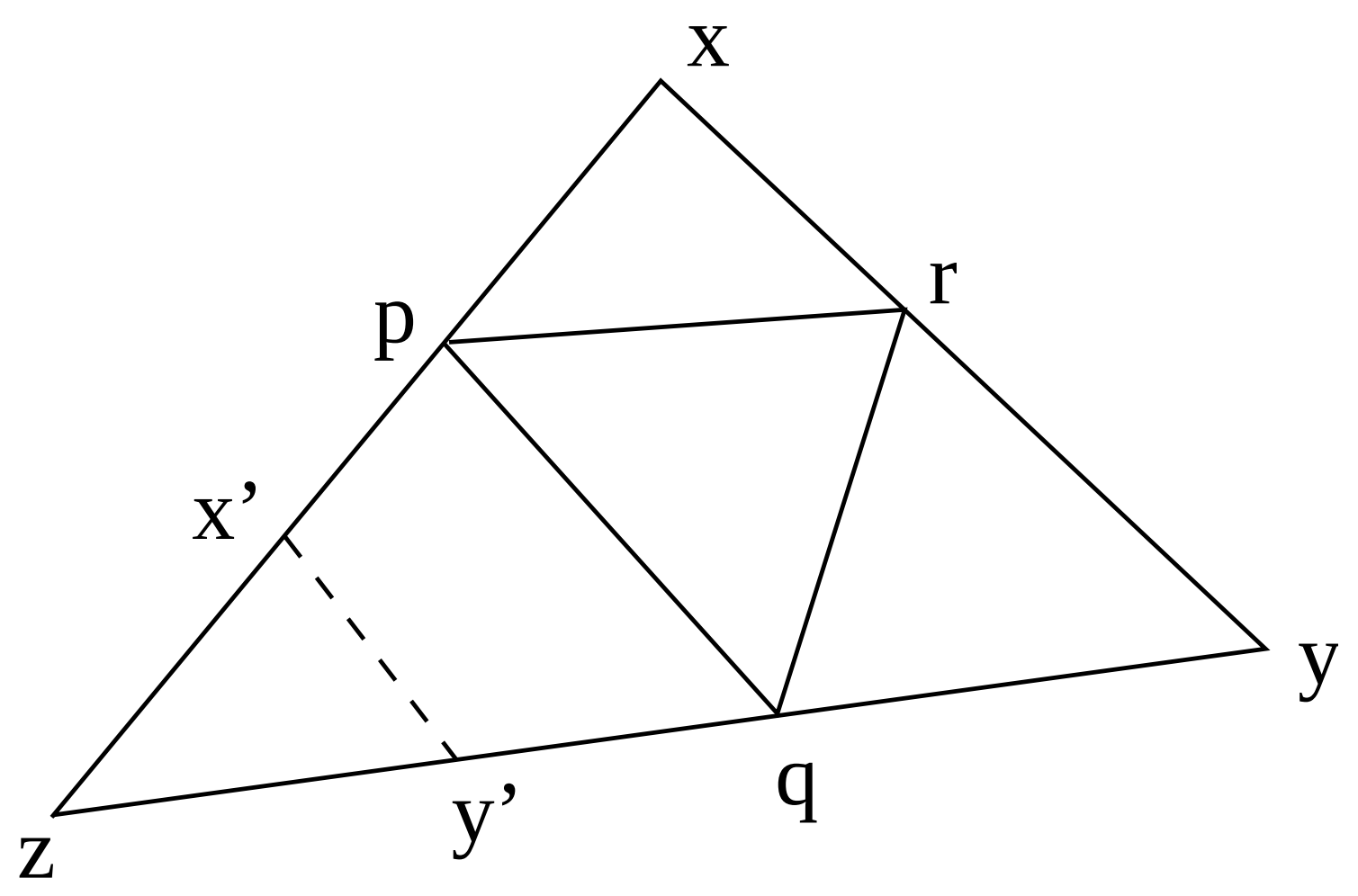

FIGURE 1. Trim triangle

We will also need another equivalent definition of hyperbolicity.

**Definition 2.4** (Gromov product)**.** Let $(X,d)$ be a metric space and suppose $x,y,z \in X$. We set

$$(x,y)_z := \frac{1}{2}[d(z,x)+d(z,y)-d(x,y)]$$

Note that $(x,y)_z = (y,x)_z$.

**Definition 2.5.** Let $(X,d)$ be a geodesic metric space and let $\Delta = \alpha_1 \cup \alpha_2 \cup \alpha_3$ is a geodesic triangle in $X$, where $\alpha_1 = [z,x]$, $\alpha_2 = [z,y]$ and $\alpha_3 = [x,y]$.

Note that by definition of Gromov product $d(x,y) = (z,y)_x + (x,z)_y$, $d(x,z) = (y,z)_x + (x,y)_z$ and $d(y,z) = (x,z)_y + (x,y)_z$. Thus there exist uniquely defined points $p \in \alpha_1$, $q \in \alpha_2$ and $r \in \alpha_3$ such that:

$$d(z,p) = d(z,q) = (x,y)_z, d(x,p) = d(x,r) = (y,z)_x, d(y,q) = d(y,r) = (x,z)_y$$

We will call $(p,q,r)$ the *inscribed triple* of $\Delta$.

**Definition 2.6** (Trim triangle)**.** Let $(X,d)$ be a geodesic metric space and let $\delta \geq 0$. Let $\Delta = \alpha_1 \cup \alpha_2 \cup \alpha_3$ be a geodesic triangle in $X$, where $\alpha_1 = [z,x]$, $\alpha_2 = [z,y]$ and $\alpha_3 = [x,y]$. We say that $\Delta$ is $\delta$*-trim* if the following holds.

Let $(p,q,r)$ be the inscribed triple of $\Delta$, where $p \in \alpha_1, q \in \alpha_2, r \in \alpha_3$, as shown in Figure 1.

Then:

1. for any points $x' \in \alpha_1, y' \in \alpha_2$ with $d(z,x') = d(z,y') \leq (x,y)_z$ we have $d(x',y') \leq \delta$;
2. for any points $z' \in \alpha_1, y' \in \alpha_3$ with $d(x,z') = d(x,y') \leq (z,y)_x$ we have $d(z',y') \leq \delta$;
3. for any points $z' \in \alpha_2, x' \in \alpha_3$ with $d(y,x') = d(y,z') \leq (x,z)_y$ we have $d(x',z') \leq \delta$.

The following statement is well-known [1]:

**Theorem 2.7.** *Let $(X,d)$ be a geodesic metric space.*

*Then $(X,d)$ is hyperbolic if and only if for some $\delta \geq 0$ all geodesic triangles in $X$ are $\delta$-trim.*

Till the end of this section let $(X,d)$ be a geodesic metric space with $\delta$-trim geodesic triangles.

The following lemma immediately follows from Definition 2.5.

**Lemma 2.8.** *Let $\Delta = [z,x] \cup [z,y] \cup [x,y]$ is a geodesic triangle in $X$ and let $a$ be an arbitrary point in the side $[x,y]$ of $\Delta$. Then either there is a point $b \in [y,z]$ such that $d(y,a) = d(y,b)$, $d(a,b) \leq \delta$ or there is a point $c \in [z,x]$ such that $d(x,a) = d(x,c)$ and $d(a,c) \leq \delta$.*

**Definition 2.9.** Suppose $(X, d)$ is a geodesic metric space. We will say that a path $p$ in $X$ is a *near geodesic* if $p = p_1 p' p_2$ such that $p'$, $p_1$, $p_2$ are geodesic segments and such that $0 \le l(p_i) \le 1$ for $i = 1, 2$.

## 3. Quasiconvex subgroups of hyperbolic groups

The detailed background information on quasiconvex subgroups of hyperbolic groups can be found in [53, 22, 37, 41, 29, 25] and other sources. We will assume some familiarity of the reader with this material.

Let $G$ be a hyperbolic group with a marked finite generating set $\pi : A \to G$. Let $X = \Gamma(G, A)$ be the Cayley graph of $G$ with respect to $A$. We will denote the word-metric corresponding to $A$ on $X$ by $d_X$. Also, for $g \in G$ we will denote $|g|_X := d_X(1, g)$. Let $\delta \ge 10$ be an integer such that $(X, d_X)$ has $\delta$-trim geodesic triangles. Let $H \le G$ be a quasiconvex subgroup. Thus $H$ is finitely generated, finitely presented and word-hyperbolic. Let $E > 0$ be an integer such that $H$ is an $E$-quasiconvex subset of $X$. For an $A$-word $w$ we will denote by $l(w)$ the length of $w$ and by $\overline{w}$ the element of $G$ represented by $w$. These constants and notations will be fixed till the end of this article, unless specified otherwise. A word $w$ will be called *$d_X$-geodesic* if $l(w) = |\overline{w}|_X$. A word $w$ will be called *near geodesic* if any path in $X$ labeled by $w$ is near geodesic.

The following observation follows immediately from Definition 2.5 since $\delta \ge 10$:

**Lemma 3.1.** *Let $\alpha$ and $\beta$ be paths from a vertex $x$ to a vertex $y$ in $X$ such that $\alpha$ is a geodesic and $\beta$ is a near geodesic. Then for any point $p \in \alpha$ there is a vertex $q \in \beta$ with $d_X(p, q) \le 3\delta$. Similarly, for any point $p \in \beta$ there is a vertex $q \in \alpha$ with $d_X(p, q) \le 3\delta$.*

The following useful statement is proved in [3]:

**Lemma 3.2.** *There exists an integer constant $K = K(G, H, A) > 0$ with the following properties. Suppose $g \in G$ is shortest with respect to $d_X$ in the coset class $Hg$. Let $h \in H$ be an arbitrary element. Then:*

1. *$|hg|_X \ge |h|_X + |g|_X - K$;*
2. *the path $[1, h] \cup h[1, g]$ is $K$-Hausdorff close to $[1, hg]$.*

Let $Y = \Gamma(G/H, A)$. Thus $Y$ is a connected graph with the induced simplicial metric, which we will denote by $d_Y$.

**Lemma 3.3.** *There exists an integer constant $K_1 = K_1(H) > 0$ with the following property.*

*Let $g, f \in G$ be shortest with respect to $d_A$ elements of $Hg$ and $Hf$ respectfully. Let $w$ be the label of a near geodesic path in $Y$ from $Hg$ to $Hf$ and let $h \in H$ be such that $hf = g\overline{w}$. Then $|h|_X \le K_1$.*

*Proof.* Let $K = K(H)$ be the constant provided by Lemma 3.2.

By the choice of $w$, the word $w$ is the label of a near-geodesic path in $X$. Let $u$ be an $A$-geodesic word representing $\overline{w} \in G$.

Consider a geodesic triangle $\Delta$ is $X$ with vertices $1, hf, g$ and with the sides $\alpha = [1, g]$, $\beta = [g, hf]$ and $\gamma = [1, hf]$. Let $u$ be the label of $\beta$, so that $\overline{u} = \overline{w} \in G$. Let $p, q, r$ be the inscribed triple of $\Delta$. Thus $p \in \alpha, q \in \beta$ and $r \in \gamma$ with $d(p, q), d(p, r), d(q, r) \le \delta$. Consider also the geodesic paths $[1, h]$ and $[h, hf] = h[1, f]$ in $X$. By Lemma 3.2 there is a point $z \in \gamma$ such that $d(h, z) \le K$. Note that $|d(1, h) - d(1, z)| \le K$. Let $\zeta$ be the path from $g$ to $fh$ labeled by $w$ in $X$.

There are two cases to consider.

**Case 1.** Suppose $d(1, z) \le d(1, r)$, as shown in Figure 2.

Let $z' \in \alpha$ be such that $d(1, z) = d(1, z')$. Then $d(z, z') \le \delta$ since $\Delta$ is $\delta$-trim. Therefore

$$d(h, g) \le K + \delta + |g|_X - d(1, z') \le K + \delta + |g|_X - |h|_X + K = |g|_X - |h|_X + 2K + \delta.$$

Since $g$ is shortest in $Hg$, we have $d(h, g) \ge |g|_X$. Therefore

$$|g|_X \le d(h, g) \le |g|_X - |h|_X + 2K + \delta$$

and hence $|h|_X \le 2K + \delta$.

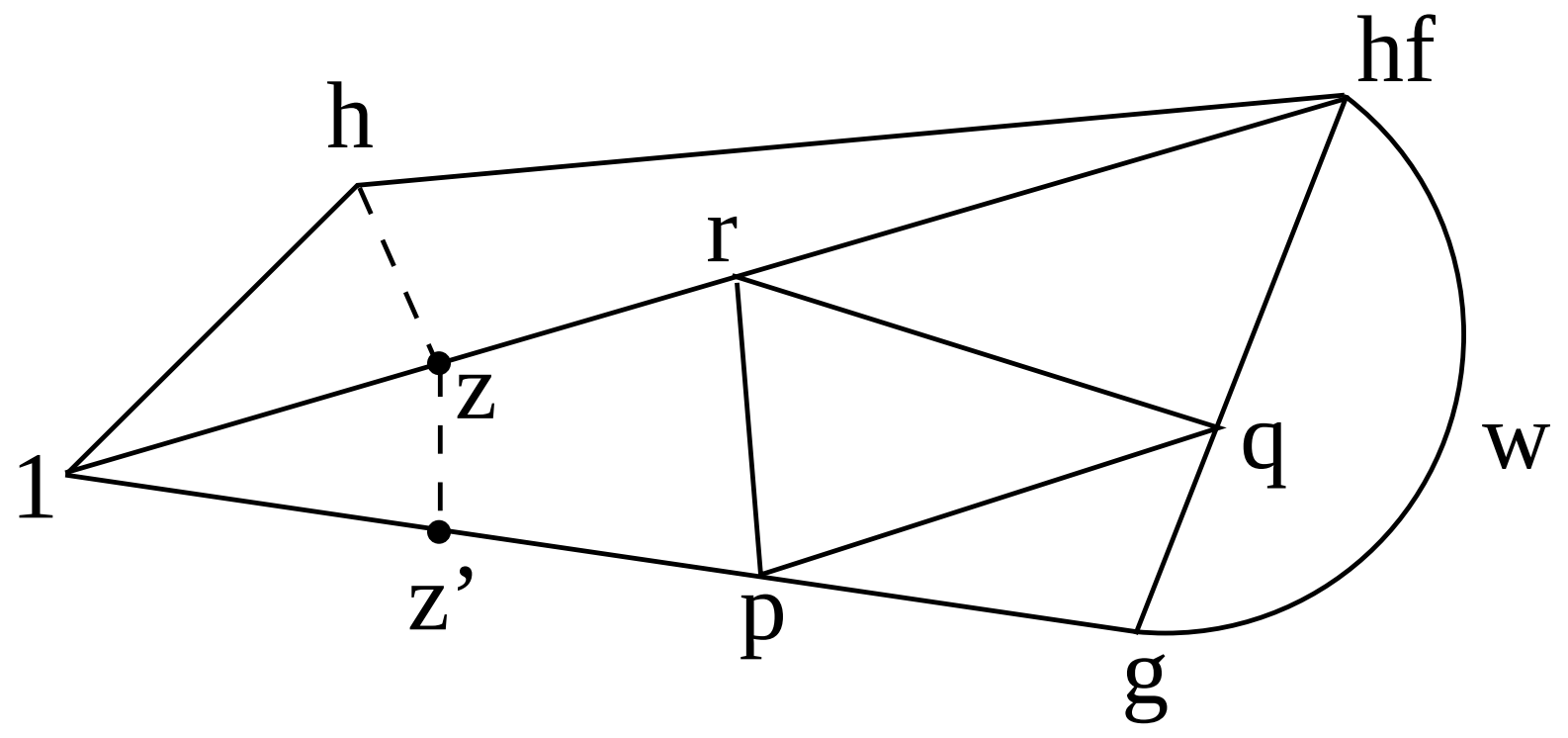

Figure 2. The case $d(1, z) \leq d(1, r)$

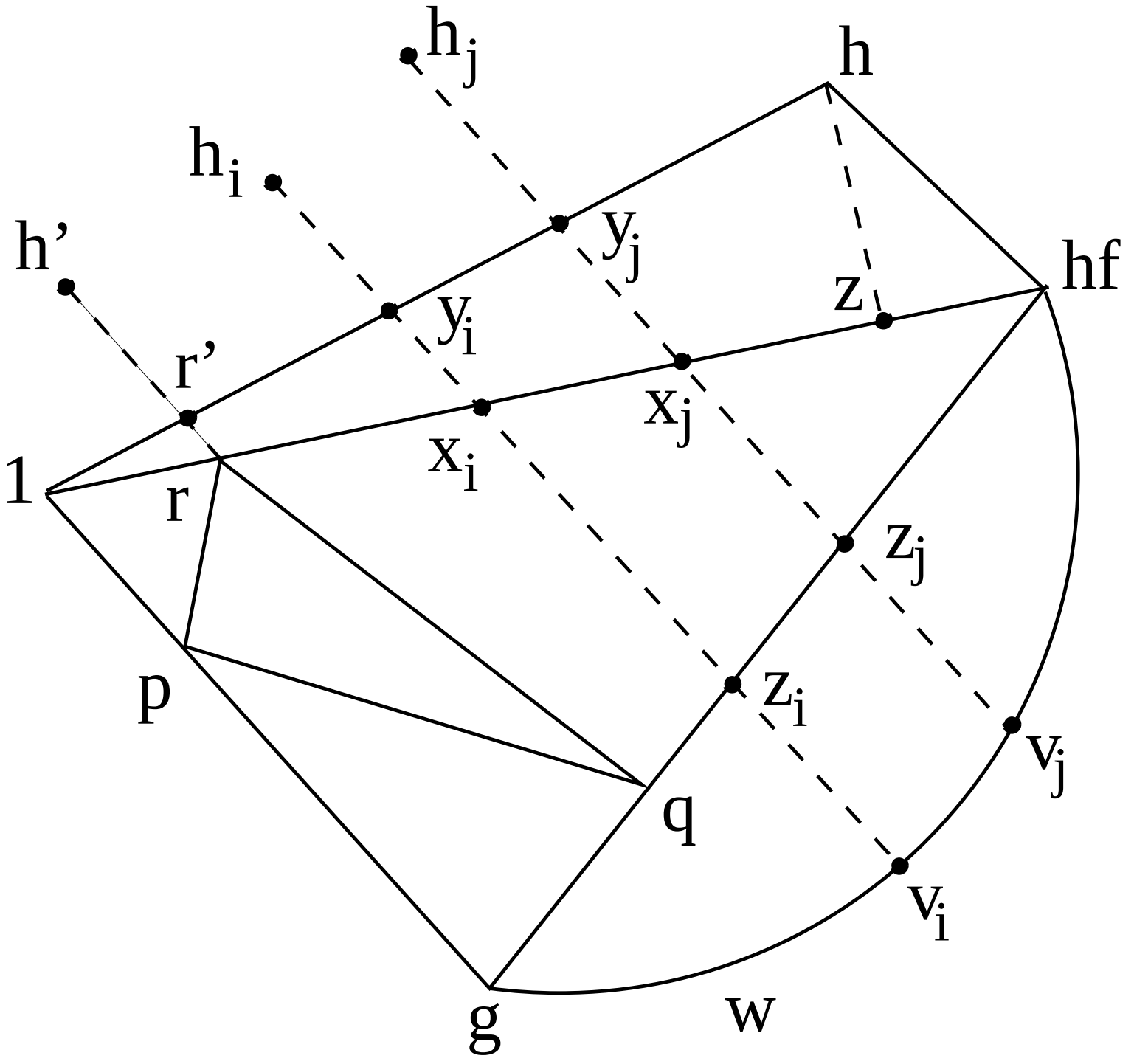

Figure 3. The case $d(1, z) > d(1, r)$

**Case 2.** Suppose $d(1, z) \geq d(1, r)$, as shown in Figure 3.

Since $d(h, z) \leq K$, there exists a point $r' \in [1, h]$ with $d(r, r') \leq K + \delta$. Since $H$ is $E$-quasiconvex in $X$, there is $h' \in H$ such that $d(r', h') \leq E$. Recall that $d(1, r) = d(1, p)$ by construction. Since $g$ is shortest in $Hg$, we have $d(h', g) \geq |g|_X$.

Therefore

$$|g|_X \leq d(h', g) \leq E + \delta + d(p, g) = E + \delta + |g|_X - d(1, r) \text{ and hence } d(1, r) \leq E + \delta.$$

Let $N$ be the number of elements in $G$ of length at most $E+K+5\delta$. Suppose that $d(r, z) \geq N(8\delta+10)$. Let $x_0 = r, x_1, \ldots, x_N$ be the points on $[r, z] \subseteq \gamma$ such that $d(r, x_i) = i(8\delta + 10)$. For each $i = 0, \ldots, N$ there is a point $y_i \in [1, h]$ such that $d(x_i, y_i) \leq K + \delta$. Since $H$ is $E$-quasiconvex in $X$, for every $i$ there is $h_i \in H$ with $d(y_i, h_i) \leq E$. (Note that we can choose $h_0 = h'$.) For every $i = 0, \dots, N$ let $z_i \in \beta$ be such that $d(hf, x_i) = d(hf, z_i)$. Then $d(x_i, z_i) \leq \delta$ since $\Delta$ is $\delta$-trim. Since $\zeta$ is a near-geodesic in $X$, for every $i$ there is a vertex $v_i$ on $\zeta$ with $d(z_i, v_i) \leq 3\delta$. Thus $d(h_i, v_i) \leq E + K + 5\delta$. By the choice of $N$ there exist $0 \leq i < j \leq N$ such that $v_i^{-1}h_i = v_j^{-1}h_j = b \in G$. Therefore $b(h_i^{-1}h_j)b^{-1} = v_i^{-1}v_j$.

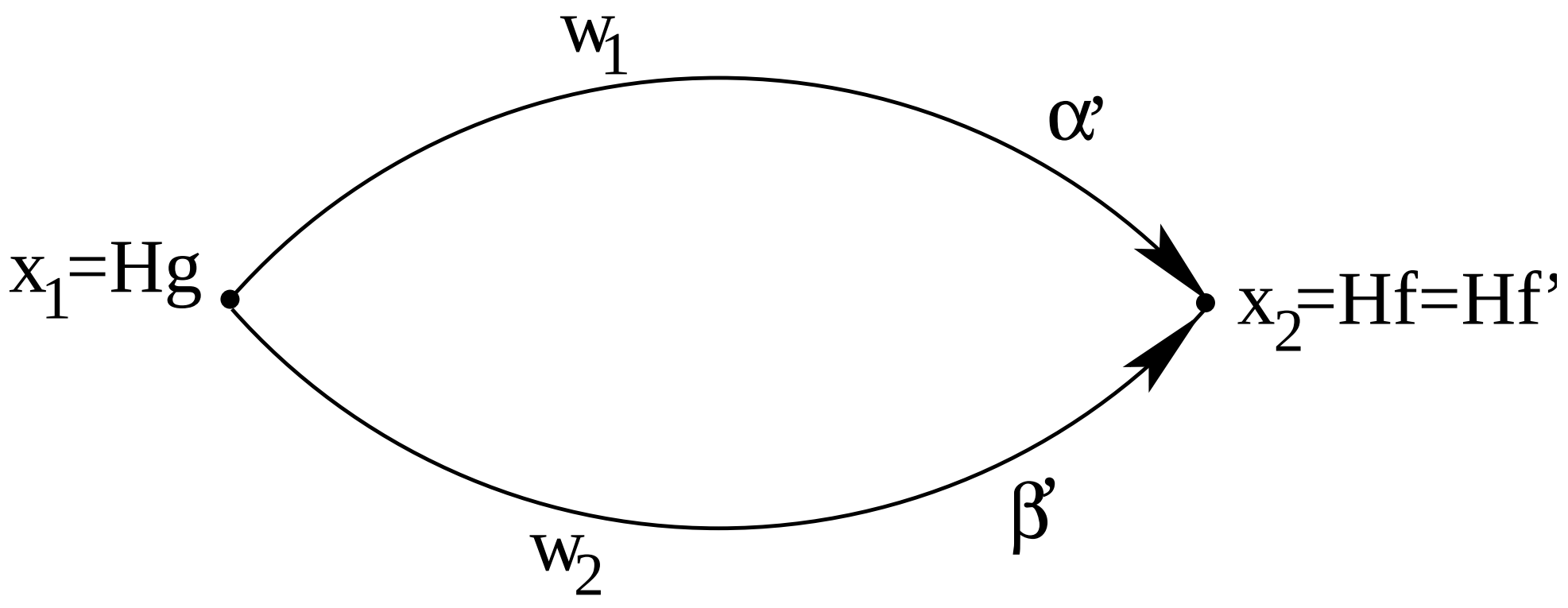

FIGURE 4. A bigon $\Theta$ in $Y$

Note that $|b|_X \leq E + K + 5\delta$ and $d(v_i, v_j) \geq 8\delta + 10 - 4\delta - 4\delta = 10$. Let $w = w_1w_2w_3$, where $w_1$ is the label of the segment of $\zeta$ from $g$ to $v_i$, where $w_2$ is the label of the segment of $\zeta$ from $v_i$ to $v_j$ and where $w_3$ is the label of the segment of $\zeta$ from $v_j$ to $hf$. Since $b(h_i^{-1}h_j)b^{-1} = v_i^{-1}v_j = \overline{w_2}$, $g\overline{w_1}b = h_i$, $g\overline{w_1w_2}b = h_j$, $bh_j^{-1}hf = \overline{w_3}$ and $v_i^{-1}h_i = v_j^{-1}h_j = b$, we have

$$h_i(h_j^{-1}h)f = g\overline{w_1w_3}$$

However $l(w_2) \geq d(v_i, v_j) \geq 10$ and hence $l(w_1w_3) \leq l(w) - 10$, contradicting the fact that $w$ is the label of a near-geodesic in $Y$ from $Hg$ to $Hf$.

Thus $d(r, z) < N(8\delta + 10)$. Since $d(1, r) \leq E + \delta$, this implies

$$|h|_X = d(1, h) \leq d(1, r) + d(r, z) + d(z, h) \leq E + \delta + N(8\delta + 10) + K.$$

Hence the statement of Lemma 3.3 holds with $K_1 := E + \delta + N(8\delta + 10) + 2K$. □

Our next goal will be to prove:

**Theorem 3.4.** *There is $\delta' \geq 0$ such that the space $(Y, d_Y)$ has $\delta'$-thin geodesic bigons.*

By Theorem 2.2 the above statement immediately implies Theorem 1.2 from the introduction.

## 4. Proof of Theorem 3.4

Since $Y$ is a graph with simplicial metric, by Theorem 2.2 it suffices to show that bigons with near-geodesic sides and vertices of $Y$ as endpoints are uniformly thin in $Y$.

We shall fix a near-geodesic bigon $\Theta$ in $Y$ with vertices $x_1, x_2$ and geodesic sides $\alpha'$ and $\beta'$ joining $x_1$ to $x_2$, as shown in Figure 4. We can write $x_1 = Hg$ and $x_2 = Hf$, where each $g, f \in G$ are shortest with respect to $d_X$ elements of the $H$-cosets $Hg$ and $Hf$ accordingly. Let the $A$-words $w, u$ be the labels of $\alpha'$ and $\beta$ accordingly. Thus $w$ and $u$ are near geodesic words with respect to $d_X$ and $Hg\overline{w} = Hg\overline{u} = Hf$. Hence there are $h_1, h_2 \in H$ such that $h_1g\overline{w} = f = h_2g\overline{u}$. By Lemma 3.3 we have $|h_i|_X \leq K_1$ for $i = 1, 2$.

Note also that $f' := g\overline{w} = hg\overline{u}$, where $h = h_1^{-1}h_2$. Thus $Hf' = Hf = x_2$ and $|h|_X \leq 2K_1$.

Consider a geodesic quadrilateral $Q$ in $X$ with vertices $1, g, hg, f'$ and geodesic sides $[1, g]$, $[1, hg]$, $\alpha = [g, f']$ and $\beta = [hg, f']$. We also draw a geodesic $\sigma = [g, hg]$, which is a diagonal of $Q$. In addition, consider the geodesics $[1, h]$ and $h[1, g] = [h, hg]$. We will work with geodesic triangles $\Delta_1 := [1, g] \cup [g, hg] \cup [1, hg]$ and $\Delta_2 := \sigma \cup \alpha \cup \beta$. Recall that $w_1, w_2$ are near-geodesic words in $A$. We will draw a path $W_1$ from $g_1$ to $f'$ with label $w_1$ and a path $W_2$ with label $w_2$ from $hg$ to $f'$. This situation is shown in Figure 5.

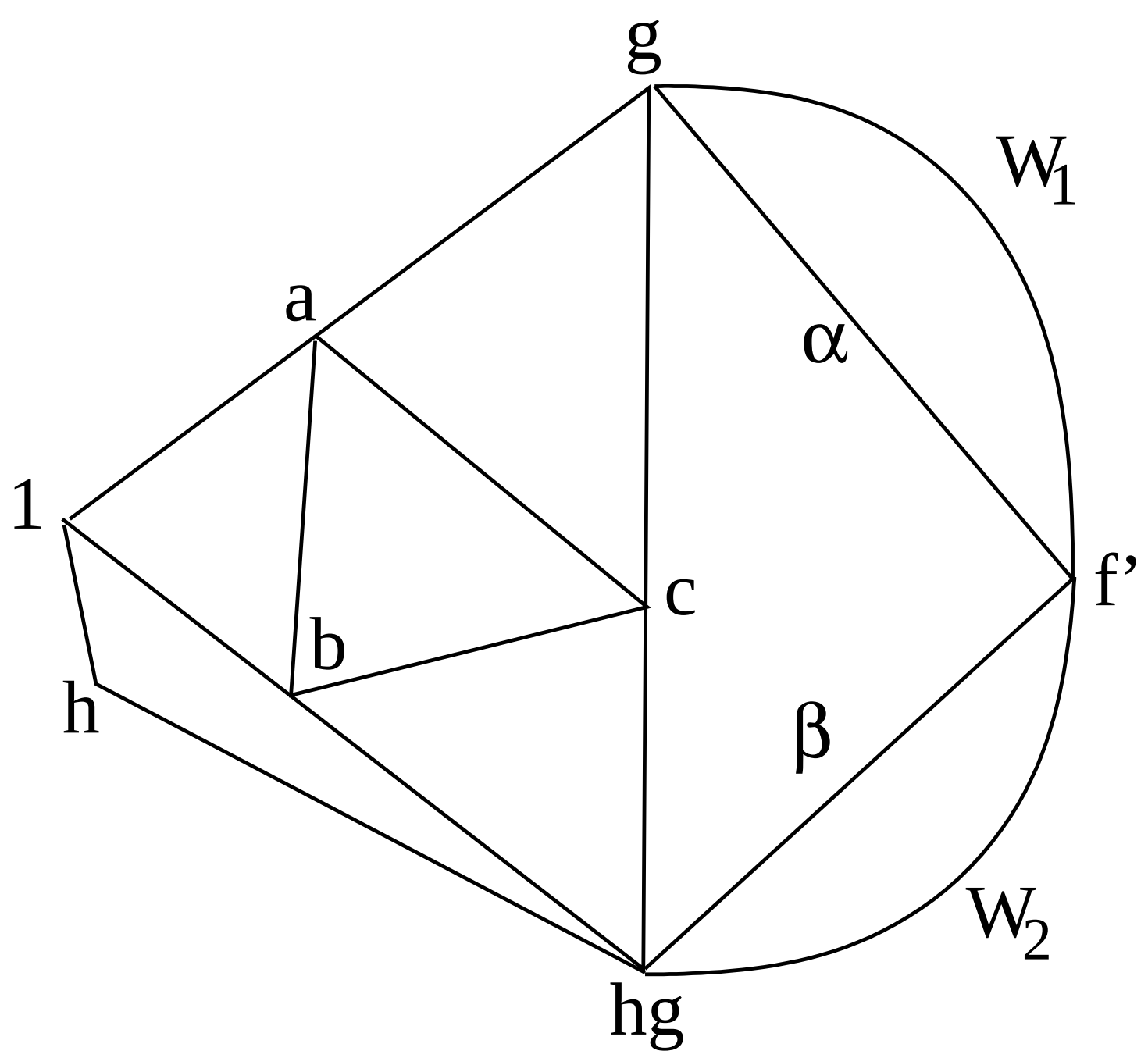

FIGURE 5. The picture in $X = \Gamma(G, A)$ corresponding to the bigon $\Theta$

Let $a, b, c$ be the inscribed triple of the triangle $\Delta_1$, where $c \in \sigma$, $a \in [1, g]$, $b \in [1, hg]$.

**Lemma 4.1.** *Let $K_1 = K_1(H) > 0$ be the constant provided by Lemma 3.3. Then $|d(g, c) - d(hg, c)| \leq 2K_1$, so that $c$ is at most $2K_1$-away from the midpoint of $\sigma = [g, hg]$.*

*Proof.* By construction $d(1, a) = d(1, b)$. Since $|h|_X \leq 2K_1$, we have $|d(1, g) - d(1, hg)| \leq 2K_1$. Hence $d(c, g) = d(a, g) = |g|_X - d(1, a)$ and $d(hg, c) = d(hg, b) = |hg|_X - d(1, b) = |hg|_X - d(1, a)$. Therefore $|d(c, g) - d(c, hg)| \leq 2K_1$, as required. □

**Lemma 4.2.** *There is a constant $K_2 = K_2(H) > 0$ with the following property. Suppose $\zeta$ is a subsegment of $\sigma = [g, hg]$ such that $c$ is the midpoint of $\zeta$ and that $\zeta$ is contained in the $\delta$-neighborhood of either $\alpha$ or $\beta$. Then $l(\zeta) \leq K_2$.*

*Proof.* Let $K_1 > 0$ be the constant provided by Lemma 3.3.

Suppose $\zeta$ be a subsegment $\sigma = [g_1, hg_1]$ such that $c$ is the midpoint of $\zeta$ and that

$$l(\zeta) \geq K_2 := 18\delta + 10K_1 + 11.$$

We will show that $\zeta$ cannot be contained in the $\delta$-neighborhood of $\beta$. The argument for $\alpha$ is completely analogous. Indeed, suppose on the contrary, that $\zeta$ is contained in the $\delta$-neighborhood of one of $\beta$. Let $z_1$ and $z_2$ be the initial and the terminal points of $\zeta$ accordingly, as shown in Figure 6.

Let $y_1 \in [g, a]$ and $y_2 \in [hg, b]$ be such that $d(g, y_1) = d(g, z_1)$ and $d(hg, y_2) = d(hg, z_2)$. Since $\Delta_1$ is $\delta$-trim, $d(y_i, z_i) \leq \delta$ for $i = 1, 2$. Moreover, as we have seen before $|d(c, g) - d(c, hg)| \leq 2K_1$. Since $d(c, z_1) = d(c, z_2)$ and $d(1, g) = d(h, hg)$, this implies that $|d(g, y_1) - d(hg, y_2)| \leq 2K_1$. Moreover, $||g|_X - |hg|_X| \leq 2K_1$ implies $|d(1, y_1) - d(1, y_2)| \leq 4K_1$.

Let $y_1' \in [h, hg] = h[1, g]$ be such that $d(1, y_1) = d(h, y_1')$. Then there is a point $t \in [1, hg]$ with $d(y_1', t) \leq 2K_1 + \delta$. Hence $|d(1, t) - d(h, y_1')| \leq 4K_1 + \delta$. It now follows from $|d(1, y_1) - d(1, y_2)| \leq 4K_1$ that $|d(1, t) - d(1, y_2)| \leq 8K_1 + \delta$. Hence $d(t, y_2) \leq 8K_1 + \delta$. Since $w_2$ is a near-geodesic in $X$ and $\zeta$ is contained in the $\delta$-neighborhood of $\beta$, there are points $p_1, p_2 \in W_2$ such that $d(z_i, p_i) \leq 4\delta$ for $i = 1, 2$. Thus $d(p_1, y_1) \leq 4\delta + \delta = 5\delta$ and

$$d(y_1', p_2) \le 2K_1 + \delta + 8K_1 + \delta + \delta + 4\delta = 10K_1 + 7\delta.$$

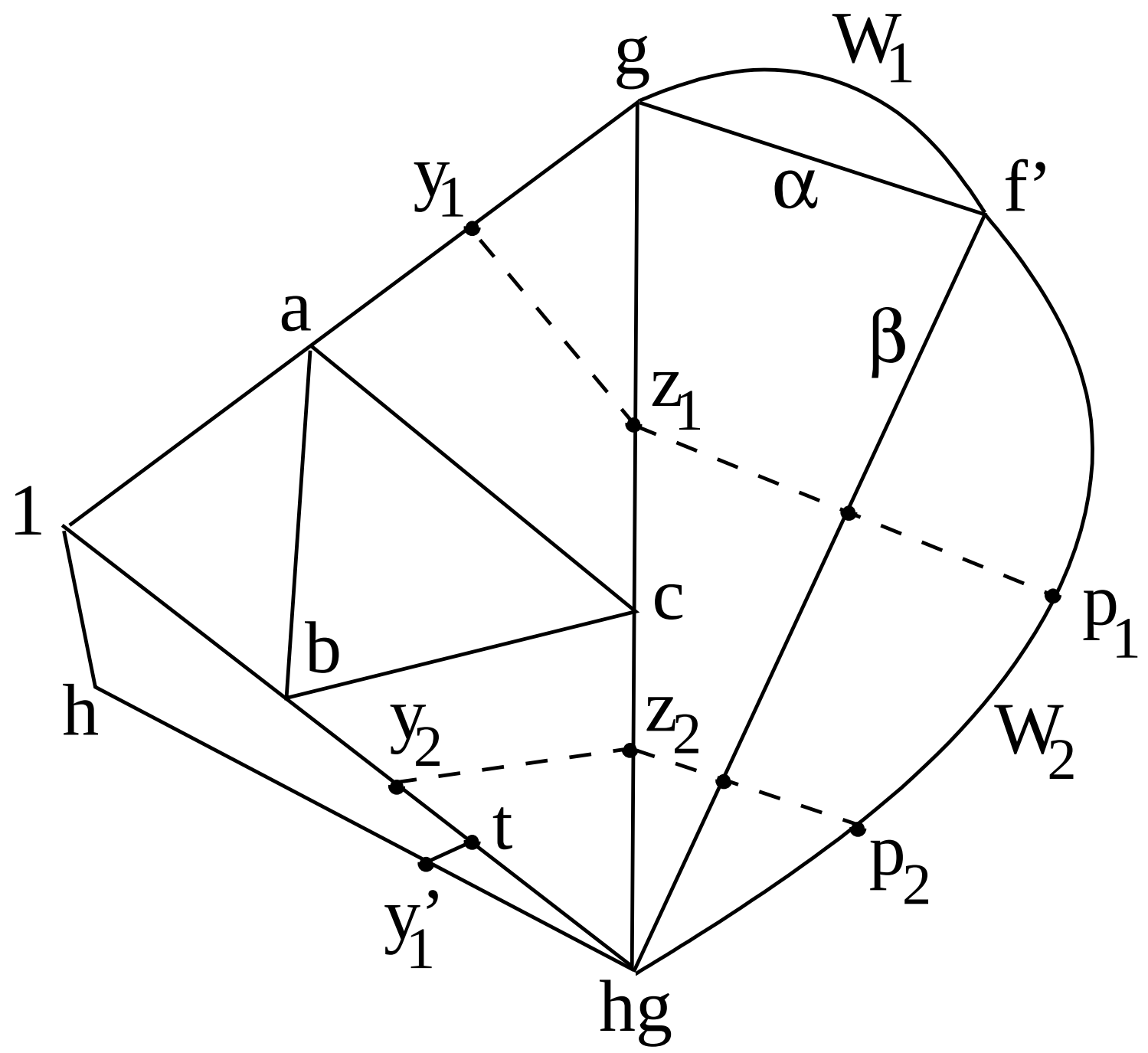

FIGURE 6. Existence of the constant $K_2$

Recall that $w_2$ is the label of $\beta$. We write $w_2$ as $w_2 = v_1 v v_2$, where $v_1$ is the label of the segment of $\beta$ from $ng$ to $p_2$, $v$ is the label of the segment of $\beta$ from $p_2$ to $p_1$ and where $v_2$ is the label the segment of $\beta$ from $p_1$ to $f'$. Hence $l(v) = d(p_1, p_2) \ge d(z_1, z_2) - 6\delta \ge K_2 - 6\delta$.

Let $u_1$ be the label of a $d_X$-geodesic path from $p_2$ to $y_1'$, so that $l(u_1) \le 7\delta + 10K_1$. Also, let $u_2$ be the label of a $d_X$-geodesic word from $y_1$ to $p_1$, so that $l(u_2) \le 5\delta$. Put $w_2' = v_1 u_1 u_2 v_2$. By construction $w_2'$ is the label of a path in $Y = \Gamma(G/H, A)$ from $x_1 = Hg$ to $x_2 = Hf = Hf'$. However $l(v) \ge K_2 - 6\delta$ and $l(u_1 u_2) \le 12\delta + 10K_1$, so that $l(u_1 u_2) < l(v) - 10$. This implies that $l(w_2') < l(w_2) - 10$, which contradicts the fact that $\beta'$ is a near geodesic in $Y$ from $x_1$ to $x_2$.

Thus $l(\zeta) \le K_2$, which completes the proof of Lemma 4.2. □

Now let $p, q, r$ be the inscribed triple for $\Delta_2$, where $p \in \sigma = [g, hg]$, $q \in \alpha = [g, f']$ and $r \in \beta = [hg, f']$.

**Lemma 4.3.** *Let $K_2 > 0$ be the constant provided by Lemma 4.2. There is an integer constant $K_3 = K_3(H) > 0$ such that either $d(g, hg) \le K_3$ or $d(p, c) \le K_2$.*

*Proof.* Put $K_3 := 2K_1 + K_2$ where $K_1$ and $K_2$ are the constants provided by Lemma 3.3 and Lemma 4.2 accordingly.

Suppose $d_X(p, c) > K_2$. We will assume that $p \in [c, g]$ as the case $p \in [c, hg]$ is symmetric. Then by Lemma 4.2 $d(gh, c) \le K_2$. Lemma 4.1 implies that $|d(g, c) - d(hg, c)| \le 2K_1$. Hence $d(g, c) \le K_2 + 2K_1$ and therefore

$$d(g, hg) = d(g, c) + d(c, hg) \le K_2 + 2K_1 + 2K_1 = K_2 + 2K_1 = K_3$$

as required. □

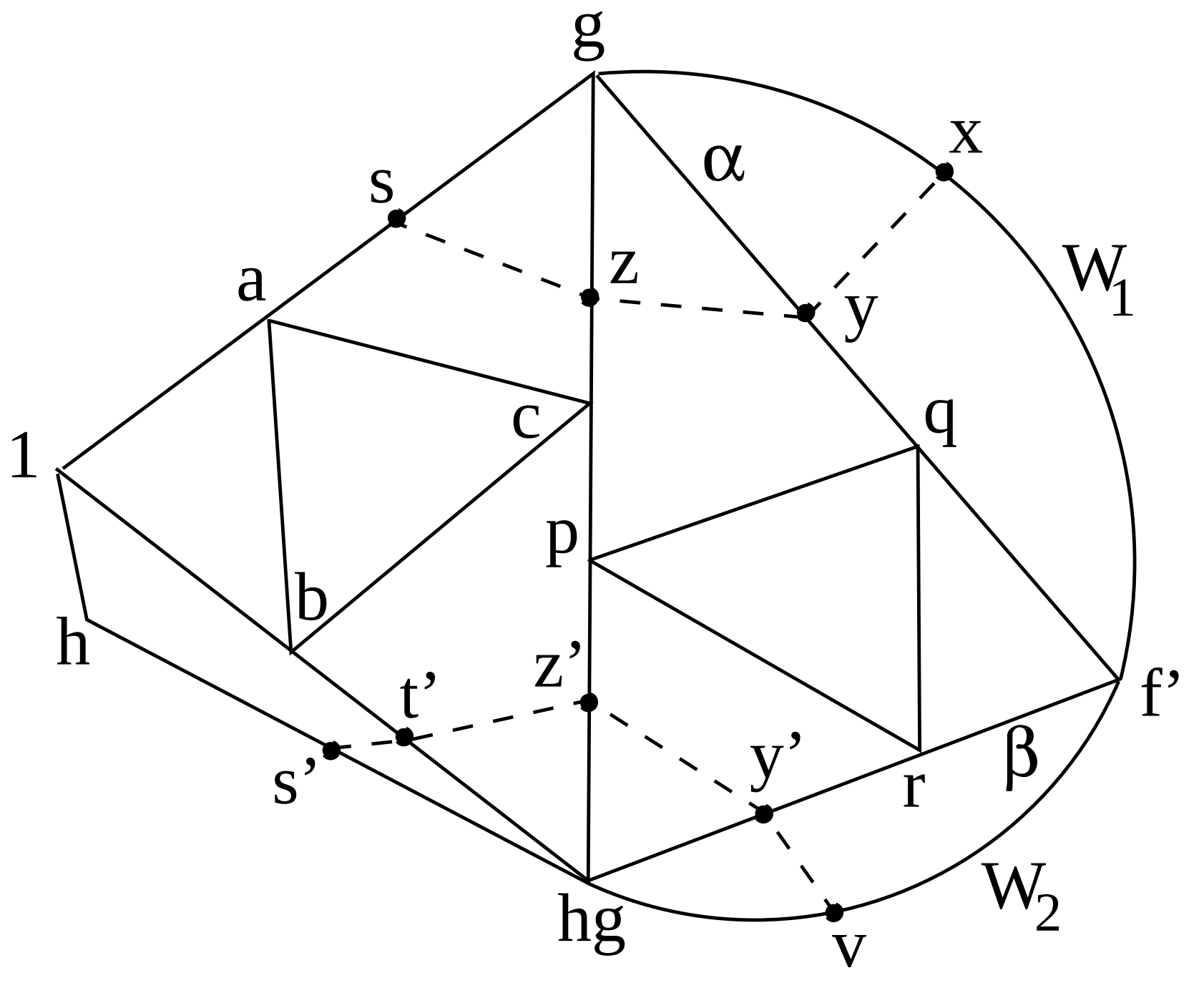

FIGURE 7. Thinness of a bigon in $Y$

*Finishing the proof of Theorem 3.4.* Suppose first that $d(g, hg) \leq K_3$, where $K_3$ is the constant provided by Lemma 4.3. Then $\alpha$ and $\beta$ are $(K_3+\delta)$-Hausdorff close. Recall that $\alpha$ is $3\delta$-close to $W_1$ and $\beta$ is $3\delta$-close to $W_2$. Hence the sides $\alpha'$, $\beta'$ of the bigon $\Theta$ in $Y$ are $(K_3+7\delta)$-Hausdorff close.

Suppose now that $d(g, hg) \geq K_3$, so that $d(p, c) \leq K_2$ by Lemma 4.3. We will assume that $p \in [hg, c]$ as the case $p \in [c, g]$ is symmetric.

Let $x' \in \alpha'$ be an arbitrary vertex. Let $x \in W_1$ be the vertex such that the segment of $W_1$ from $g$ to $x$ has the same label as the segment of $\alpha'$ from $x_1$ to $x'$. There is a vertex $y \in \alpha$ such that $d(x, y) \leq 3\delta$. Suppose first that $y \in [q, f']$. Since $\Delta_2$ is $\delta$-trim, there is a vertex $z$ on $W_2$ such that $d(x, z) \leq 4\delta$. Hence $d(x, z) \leq 7\delta$. Let $z' \in \beta'$ be such that the segment of $W_2$ from $z$ to $f'$ and the segment of $\beta'$ from $z'$ to $x_2$ have the same label. Then $d_Y(x', z') \leq 7\delta$.

Suppose now that $y \in [g, q]$, as shown in Figure 7. Recall that $d(p, c) \leq K_2$. Since $\Delta_2$ is $\delta$-trim, there is a point $z \in [c, g]$ with $d(y, z) \leq K_2 + \delta$. Let $s \in [a, g]$ be such that $d(g, s) = d(g, z)$ and so $d(z, s) \leq \delta$. Let $s' \in [h, hg] = h[1, g]$ be such that $d(h, s') = d(1, s)$. Since $d(1, h) \leq 2K_1$, there is a point $t \in [1, hg]$ with $d(s', t) \leq 2K_1 + \delta$. Hence $|d(1, t) - d(1, s)| \leq 4K_1 + \delta$. Since $d(1, a) = d(1, b)$ and $d(1, s) \geq d(1, a)$, there is a point $t' \in [b, hg]$ with $d(t, t') \leq 4K_1 + \delta$. Recall that $d(c, p) \leq K_2$. Since $\Delta_1$ is $\delta$-trim, there exists a point $z' \in [hg, p]$ such that $d(t', z') \leq K_2 + \delta$. Because $\Delta_2$ is $\delta$-trim, there is a point $y' \in [hg, r]$ with $d(z', y') \leq \delta$. Finally, there exists a vertex $v$ on $W_2$ such that $d(y', v) \leq 3\delta$.

Thus

$$d(x, s) \leq d(x, y) + d(y, z) + d(z, s) \leq 3\delta + K_2 + \delta + \delta = K_2 + 5\delta$$

and

$$\begin{aligned} d(s', v) &\leq d(s', t) + d(t, t') + d(t', z') + d(z', y') + d(y', v) \leq \\ &\leq (2K_1 + \delta) + (4K_1 + \delta) + (K_2 + \delta) + \delta + 3\delta = 6K_1 + K_2 + 7\delta. \end{aligned}$$

Let $v' \in \beta'$ be the vertex such that the initial segment of $\beta'$ from $x_1$ to $v'$ and the initial segment of $W_2$ from $hg$ to $v$ have the same label.

By the choice of $s$ and $s'$ we have

$$d(x', v') \leq K_2 + 5\delta + 6K_1 + K_2 + 7\delta = 6K_1 + 2K_2 + 12\delta.$$

Thus we have shown that $\alpha'$ is contained in the $\delta'$-neighborhood of $\beta'$, where

$$\delta' := \max\{6K_1 + 2K_2 + 12\delta, K_3 + 7\delta\}.$$

A virtually identical argument shows that $\beta'$ is contained in the $\delta''$-neighborhood of $\alpha'$, where $\delta'' > 0$ is some constant independent of the choice of the bigon $\Theta$. This completes the proof of Theorem 3.4. □

## 5. Simple random walks

**Definition 5.1.** Let $X$ be a connected graph of bounded degree. We will say that $X$ is *recurrent* if the simple random walk on $Y$ eventually returns to the base-point with probability 1. Otherwise $X$ is said to be *transient.*

We refer the reader to [9, 57] for the detailed background information about random walks on graphs and for further references on this vast subject.

**Definition 5.2.** Recall that for two metric spaces $(X, d_X)$ and $(Y, d_Y)$ a map $f : X \to Y$ is said to be *quasi-isometry* if there is $C > 0$ such that:

(i) For any $x_1, x_2 \in X$ we have

$$\frac{1}{C} d_X(x_1, x_2) - C \leq d_Y(f(x_1), f(x_2)) \leq C d_X(x_1, x_2) + C.$$

(ii) For any $y \in Y$ there is $x \in X$ with $d(y, f(x)) \leq C$.

If only condition (i) above is known to hold, the map $f$ is said to be a *quasi-isometric embedding* (even though $f$ need not be injective).

We need the following simple fact which asserts that a graph of bounded degree, which admits a quasi-isometric embedding of a binary tree, is transient.

**Proposition 5.3.** *Let $B$ be the regular binary rooted tree with simplicial metric $d_B$ (that is the root-vertex has degree two and every other vertex has degree three). Let $Y$ be a graph of bounded degree with simplicial metric $d_Y$ which admits a quasi-isometric-embedding $VB \to VY$. That is, suppose there exists a map $f : VB \to VY$ such that for some $C > 0$ we have*

$$\frac{1}{C} d_B(v, u) - C \leq d_Y(f(v), f(u)) \leq C d_B(u, v) + C$$

*for any vertices $v, u \in VB$. Then $Y$ is transient.*

*Proof.* Choose a root vertex $v_0$ in $B$ and orient the edges of $B$ so that the edges pointing away from $v_0$ are positive. For every positive edge $e = [u, v]$ of $B$ choose a geodesic segment $J_e = [f(u), f(v)]$ in $Y$. Note that $l(J_e) \leq 2C$ by the properties of $f$. Put $Z$ to be the union of all $J_e$ where $e$ varies over all positive edges of $B$. Then $Z$ is obviously a connected subgraph of $Y$ (and hence has bounded degree). Let $d_Z$ be the simplicial metric for $Z$ (which may differ from $d_Y|_Z$).

We claim that $Id : (VZ, d_Z) \to (VY, d_Y)$ is a quasi-isometry. Since part (ii) of Definition 5.2 obviously holds, we need to show that there exists $C_1 > 0$ such that for any vertices $z_1, z_2$ of $Z$ we have

$$(\dagger) \qquad \frac{1}{C_1} d_Y(z_1, z_2) - C_1 \leq d_Z(z_1, z_2) \leq C_1 d_Y(z_1, z_2) + C_1.$$

Since every vertex of $Z$ is connected in $Z$ by a path of length at most $C$ to a vertex in $f(VB)$, it suffices to establish the above inequality for $z_1, z_2 \in F(VB)$.

If $z_1 = f(u_1)$, $z_2 = f(u_2)$ then $d_Y(z_1, z_2) = d_Y(f(u_1), f(u_2)) \geq \frac{1}{C} d_B(u_1, u_2) - C$, that is $d_B(u_1, u_2) \leq C d_Y(f(u_1), f(u_2)) + C^2$. On the other hand since the distances between the images of adjacent vertices of $B$ are bounded by $2C$, we have

$$d_Z(f(u_1), f(u_2)) \leq 2C d_B(u_1, u_2).$$

Combining these inequalities, we obtain

$$d_Z(z_1, z_2) \leq 2Cd_B(u_1, u_2) \leq 2C^2 d_Y(f(u_1), f(u_2)) + C^3.$$

On the other hand, it is obvious that $d_Y(z_1, z_2) \leq d_Z(z_1, z_2)$. Hence for $z_1 = f(u_1), z_2 = f(u_2)$ the inequality (†) holds with $C_1 = \max\{1, 2C^2, C^3\}$.

Thus we have established that (†) holds and so $f : (VB, d_B) \to (VZ, d_Z)$ is a quasi-isometry. Since $B$ is a regular binary tree, it is transient. By Theorem 3.10 of [57] transience is a quasi-isometry invariant for connected graphs of bounded degree. Hence $Z$ is a transient graph. Since $Z$ is a subgraph of $Y$ and $Y$ is locally finite, Lemma 2.2 of [7] implies that $Y$ is transient as well. □

**Theorem 5.4.** *[c.f. Theorem 1.3 from the Introduction.] Let $G$ be a torsion-free non-elementary word-hyperbolic group with a marked finite generating set $A$. Let $H \leq G$ be a quasiconvex subgroup and let $Y = \Gamma(G/H, A)$ be the relative Cayley graph for $G$ relative $H$.*

*Then $Y$ is transient if and only if $[G : H] = \infty$.*

*Proof.* If $[G : H] < \infty$, then $Y$ is finite and hence clearly recurrent. Suppose now that $[G : H] = \infty$. If $H = 1$ then $G = G/H$ and the statement is obvious since by assumption $G$ is non-elementary and thus contains a free quasiconvex subgroup of rank two [23, 10]. Hence there is a quasi-isometric embedding of a 4-regular tree into $Y = \Gamma(G, A)$ and so $Y$ is transient by Proposition 5.3.

Suppose now $H \neq 1$. Let $X = \Gamma(G, A)$ be the Cayley graph of $G$ with respect to $A$. We will denote the word-metric on $X$ corresponding to $A$ by $d_X$. Also for $g \in G$ we denote $|g|_X := d_X(1, g)$. Put $Y = \Gamma(G/H, A)$. Thus $Y$ is a connected $2k$-regular graph where $k$ is the number of elements in $A$. We denote the simplicial metric on $Y$ by $d_Y$.

By a result of G.Arzhantseva [2] since $H \leq G$ is quasiconvex and has infinite index, there exists an element of infinite order $c \in G$ such that $H_1 = \langle H, c\rangle = H * \langle c\rangle$ is quasiconvex in $G$. Choose $h_0 \in H - \{1\}$ and put $a = ch_0c^{-1}$, $b = c^2h_0c^{-2}$. It is easy to see that $H_2 = \langle H, a, b\rangle = H * F(a, b) \leq H_1 \leq G$ is quasiconvex in $H_1$ and hence in $G$. Thus $H$ is quasiconvex in $H_2$, $H_1$ and $G$. This implies that the free group $F(a, b)$ with the standard metric admits a quasi-isometric embedding in $Y$. Indeed, consider the map $\xi : F(a, b) \to Y$ defined as $\xi : f \mapsto Hf$ for $f \in F(a, b)$.

For $f \in F := F(a, b)$ we will denote by $|f|_F$ the freely reduced length of $f$ in $F(a, b)$ with respect to the free basis $a, b$. We will also denote by $d_F$ the word-metric on $F(a, b)$ corresponding to the free basis $a, b$. Choose a finite generating set $Q$ for $H$ and denote by $d_Q$ the word-metric on $H$ corresponding to $Q$. For $h \in H$ we will denote $|h|_Q := d_Q(1, h)$. Put $S = \{a, b\} \cup Q$ so that $S$ is a generating set for $H_2$. Denote by $d_S$ the word metric on $H_2$ corresponding to $S$. Since $H_2$ is quasiconvex in $G$, there is a constant $\lambda > 0$ such that for any element $g \in H_2$ we have $d_X(1, g) \geq \lambda d_S(1, g)$. Moreover, since $F$ is a finitely generated subgroup of $G$, there is a constant $\lambda' > 1$ such that for any $f_1, f_2 \in F$ we have $d_X(f_1, f_2) \leq \lambda' d_F(f_1, f_2)$. (In fact we can choose $\lambda' = \max\{|a|_X, |b|_X\}$.)

Suppose $f_1, f_2 \in F$ are two arbitrary elements. Then $\xi(f_i) = Hf_i$. Let $w$ be the $A$-word which is the label of a geodesic path in $Y$ from $Hf_1$ to $Hf_2$. Denote by $g$ the element of $G$ represented by $w$. Then there is $h \in H$ such that $f_1g = hf_2$. Thus $g = f_1^{-1}hf_2 \in H_2$. Since $H_2 = H * F$, this implies

$$d_S(1, g) \geq |f_1^{-1}f_2|_F + |h|_Q \geq |f_1^{-1}f_2|_F = d_F(f_1, f_2).$$

This implies

$$d_Y(Hf_1, Hf_2) = l(w) = d_X(1, g) \geq \lambda d_S(1, g) \geq \lambda d_F(f_1, f_2).$$

On the other hand there is a path of length $|f_1^{-1}f_2|_X$ from $f_1$ to $f_2$ in $X$. Hence there is a path of length $|f_1^{-1}f_2|_X$ from $Hf_1$ to $Hf_2$ in $Y$. Therefore $d_Y(Hf_1, Hf_2) \leq |f_1^{-1}f_2|_X = d_X(f_1, f_2) \leq \lambda' d_F(f_1, f_2)$.

Thus

$$\lambda d_F(f_1, f_2) \leq d_Y(Hf_1, Hf_2) \leq \lambda' d_F(f_1, f_2).$$

Hence $\xi$ is a quasi-isometric embedding of $(F, d_F)$ into $(Y, d_Y)$. The Cayley graph of $F$ is a 4-regular tree which contains a regular binary tree as an isometrically embedded subgraph. Hence by Proposition 5.3 the graph $Y$ is transient. □

Department of Mathematics, University of Illinois at Urbana-Champaign, 1409 West Green Street, Urbana, IL 61801, USA
*E-mail address*: kapovich@math.uiuc.edu